\documentclass[twoside,11pt]{book}
\usepackage{amsmath}
\usepackage{epsfig}
\usepackage{psfig}
\usepackage{latexsym}

\def\Real{{I\!\!R}}

\newtheorem{theorem}{Theorem}[section]

\newtheorem{definition}{Definition}[section]

\begin{document}

\title{The Local Convergence of the   Extended Kalman Filter}
\author{ Arthur J. Krener
\thanks{Research supported
in part by NSF DMS-0204390 and  AFOSR F49620-01-1-0202.}}

\maketitle

\chapter{The Convergence of the   Extended Kalman Filter}
\section*{Abstract}
We demonstrate that the extended Kalman filter converges locally for a broad 
class of nonlinear systems.  If the initial estimation error of the filter is not too large then 
the error goes to zero exponentially as time goes to infinity. To demonstrate this, 
we require that the system be  $C^2$ and 
uniformly observable with bounded second partial derivatives.

\section{Introduction}
\setcounter{equation}{0}
The extended Kalman filter is a widely used method for estimating the state
$x(t)\in \Real^n$ of 
a partially observed nonlinear dynamical system,
\begin{equation}\begin{array}{lllllllllllllllllllllllllllllllll} \label{nl}
\dot{x}&=& f(x,u)\\
y&=& h(x,u)\\
x(0)&=& x^0
\end{array}\end{equation}
from the past controls and  observations $u(s)\in { U}\subset\Real^m,\ y(s)\in \Real^p, 0\le s\le t$ and some
information about
the initial condition $x^0$.
The functions $f, h$ are known and assumed to be  $C^2$.

An extended Kalman filter is derived by replacing (\ref{nl}) by its
linear approximation around the trajectory $\hat{x}(t)$ and adding standard
white Gaussian  driving 
noise $w(t)\in \Real^l$ and independent, standard white Gaussian  observation 
noise $v(t)\in \Real^l$, 
\begin{equation}\begin{array}{lllllllllllllllllllllllllllllllll}  \label{lin}
\dot{z}&=&f(\hat{x}(t),u(t))+A(t)z+Gw \\
y&=&h(\hat{x}(t),u(t))+C(t)z+v\\
z(0)&=&z^0
\end{array}\end{equation}
where $G$ is a $n\times l$ matrix chosen by the designer,
\begin{equation}\begin{array}{lllllllllllllllllllllllllllllllll}
A(t)={\displaystyle \frac{\partial f}{\partial x}(\hat{x}(t),u(t))},&&
\hspace{0.25in} 
 C(t)={\displaystyle \frac{\partial h}{\partial
x}(\hat{x}(t),u(t))},
\end{array}\end{equation}
and $z^0$ is a Gaussian random vector independent of the noises
with mean $\hat{x}^0$ and variance $P^0$ that are
chosen by the designer. 

The Kalman filter for (\ref{lin}) is
\begin{equation}\begin{array}{lllllllllllllllllllllllllllllllll} \label{kf}
\dot{\hat{z}}(t)&=&f(\hat{x}(t),u(t))+A(t)\hat{z}(t)+P(t) C'\left(y(t)-h(\hat{x}(t),u(t))  -C(t)\hat{z}(t)\right)\\
\dot{P}(t)&=&A(t)P(t)+P(t)A'(t)+ \Gamma -P(t)C'(t)C(t)P(t)\\
\hat{z}(0)&=&\hat{x}^0\\
P(0)&=&P^0
\end{array}\end{equation}
where $\Gamma=GG'$.

The extended Kalman filter for (\ref{nl}) is given by
\begin{equation}\begin{array}{lllllllllllllllllllllllllllllllll}
\label{ekf}
\dot{\hat{x}}(t)&=&f(\hat{x}(t),u(t))
+P(t) C'(t)\left(y(t)-h(\hat{x}(t),u(t))\right)\\
\dot{P}(t)&=&A(t)P(t)+
P(t)A'(t)+\Gamma -P(t)C'(t)C(t)P(t)\\
\hat{x}(0)&=&\hat{x}^0\\
P(0)&=&P^0.
\end{array}\end{equation}
Actually there are many extended Kalman filters for (\ref{nl}), depending on the choice
of the design parameters $G, \ \hat{x}^0,\ P^0$.  We could also broaden 
the class of extended Kalman filters for (\ref{nl}) by allowing $G=G(\hat{x}(t))$ and putting a similar
coefficient in front of the observation noise in (\ref{lin}).  We chose not to do so to simplify
the discussion.  For similar reasons we omit the discussion of time varying systems.
We expect that our main theorem can be generalized to cover such systems.
For more on the derivation of the extended Kalman filter, see Gelb \cite{G74}.

Baras, Bensoussan and James \cite{BBJ88} have shown that under suitable conditions, the extended 
Kalman filter converges locally, i.e.,  if the initial error $\tilde{x}(0)=x(0)-\hat{x}(0)$ is sufficiently
small then $\tilde{x}(t)=x(t)-\hat{x}(t) \to 0$ as $t \to \infty$.  Unfortunately their conditions are difficult
to verify and may not be satisfied even by an observable linear system.   
Krener and Duarte have given a simple example where any   extended Kalman filter
fails to converge.  More on these points later. 

By modifying the techniques of \cite{BBJ88} and incorporating techniques of the 
high gain observer of Gauthier,  Hammouri and  Othman \cite{GHO92}
 we shall show that under  verifiable conditions that the extended Kalman 
filter converges locally. To state the main result we need a definition.
\begin{definition} \cite{GHO92} The system
\begin{equation}\begin{array}{lllllllllllllllllllllllllllllllll} \label{nl2}
\dot{\xi}&=& f(\xi,u)\\
y&=& h(\xi,u)
\end{array}\end{equation}
     is uniformly  observable for any input if there exist  coordinates
$$\left\{x_{ij}: i=1,\ldots,p,\, j=1,\ldots,
l_i\right\}$$
where $1\le l_1 \le \ldots\le l_p$ and $\sum l_i=n$  such that in
these  coordinates the system takes the
form
\begin{equation}\begin{array}{lllllllllllllllllllllllllllllllll} \label{sys1}
y_i&=&x_{i1}+h_i(u) \\
\dot{x}_{i1}&=&{x}_{i2}+f_{i1}(\underline{x}_1,u) \\
&\vdots&\\
\dot{x}_{ij}&=&{x}_{ij+1} +f_{ij}(\underline{x}_j,u) \\
&\vdots&\\
\dot{x}_{il_i-1}&=&{x}_{il_i}+f_{il_i-1}(\underline{x}_{l_i -1},u) \\
\dot{x}_{il_i}&=&f_{il_i }(\underline{x}_{l_i },u)
\end{array}\end{equation}
for $i=1,\ldots,p$ where $ \underline{x}_j$ is defined by
\begin{eqnarray}
\underline{x}_j&=& (x_{11},\ldots,x_{1,j\wedge l_1},x_{21},\ldots,x_{pj}).
\end{eqnarray}
Notice that in $ \underline{x}_j$ the indices range over
$i=1,\ldots,p;\  k=1,\ldots,\min\{j,l_i\}$ and the coordinates are
ordered so that
     second index moves faster
than the first.

     We also require that each $f_{ij} $ be Lipschitz continuous, there exists an $L$ such that
     for all $x,\xi\in \Real^n, u\in U$, 
\begin{equation}\begin{array}{lllllllllllllllllllllllllllllllll} \label{LipGr}
|f_{i}(\underline{x}_j,u)-f_{i}(\underline{\xi}_j,u)|&\le& L|\underline{x}_j-\underline{\xi}_j|.
\end{array}\end{equation}
 The symbol $|\cdot|$ denotes the Euclidean norm.  
\end{definition}

Let
\begin{equation}\begin{array}{lllllllllllllllllllllllllllllllll} \nonumber
\bar{A}_i&=&\left[ \begin{array}{ccccccccc} 0& 1& 0& \ldots& 0\\
0&0&1 &\ldots& 0\\
&&&\ddots&\\
0&0&0 &\ldots& 1\\
0&0&0 &\ldots& 0 \end{array}\right]^{l_i\times l_i}&&
\bar{A}&=& \left[ \begin{array}{ccccccccc} \bar{A}_1& 0& 0\\
0&\ddots &0\\
0&0 & \bar{A}_p \end{array}\right]^{n\times n}
\end{array}\end{equation}
\begin{equation}\begin{array}{lllllllllllllllllllllllllllllllll}
\bar{C}_i&=&\left[ \begin{array}{ccccccccc} 1&0&0&\ldots&0 \end{array}\right]^{1\times l_i}
&&\bar{C}&=& \left[ \begin{array}{ccccccccc} \bar{C}_1& 0& 0\\
0&\ddots &0\\
0&0 & \bar{C}_p \end{array}\right]^{p\times n} \nonumber
\end{array}\end{equation}
\begin{equation}\begin{array}{lllllllllllllllllllllllllllllllll}
\bar{f}_i(x,u)&=&\left[ \begin{array}{ccccccccc} f_{i1}(\underline{x}_1,u)\\ \vdots \\
f_{il_i}(\underline{x}_{l_i},u)\end{array}\right]^{l_i\times
1}&&
\bar{f}(x,u)&=& \left[ \begin{array}{ccccccccc} \bar{f}_1(x,u)\\ \vdots\\ \bar{f}_p(x,u) \end{array}\right]^{n\times 1}
\nonumber\\  \\
&&&&\bar{h}(u)&=&\left[ \begin{array}{ccccccccc} h_1(u)\\ \vdots \\ h_p(u) \end{array}\right]^{p\times 1}
\end{array}\end{equation}
then (\ref{sys1}) becomes
\begin{equation}\begin{array}{lllllllllllllllllllllllllllllllll} \label{nl3}
\dot{x}&=& \bar{A}x+\bar{f}(x,u)\\
y&=&\bar{C}x+\bar{h}(u)
\end{array}\end{equation}
A system such as (\ref{sys1}) or, equivalently (\ref{nl3}), is said to be in observable form \cite{K86}.

We shall also require that the second derivative of $\bar{f}$ is bounded, i.e., for any
$x,\xi\in \Real^n, u\in U$,
\begin{eqnarray} \label{secpar}
\left| \frac{\partial^2 \bar{f}}{\partial x_i \partial x_j}(x,u)\xi_i\xi_j\right|\le L |\xi|^2.
\end{eqnarray}
On the left we employ the convention of summing on repeated indices. 

\begin{theorem} {\bf (Main Theorem)} \label{mt}
Suppose
\begin{itemize}
\item the system (\ref{nl}) is uniformly observable for any input and so without
loss of generality we can assume that is in the
form (\ref{nl3}) and satisfies the Lipschitz conditions (\ref{LipGr}), 
\item  the second derivative of $\bar{f}$ is bounded (\ref{secpar}),
\item $x(t), y(t)$ are any state and output trajectories generated by
(\ref{nl3}), 
\item $G$ has been chosen to be invertible,
\item $\hat{x}(t)$ and $P(t)$ are a solution of the extended Kalman filter 
(\ref{ekf}) where $P(0)$ is positive definite and $\tilde{x}(0)=x(0)-\hat{x}(0)$ is sufficiently small,
\end{itemize}
\vspace{0.05in}
Then  $|x(t)-\hat{x}(t)| \to 0$ exponentially as $t \to \infty$.  
\vspace{0.1in}
\end{theorem}

\section{Proof of the Main Theorem}
\setcounter{equation}{0}
We extend the method of proof of \cite{BBJ88}. 
Since the system is in observable form
\begin{equation}\begin{array}{lllllllllllllllllllllllllllllllll}
A(t)&=&\bar{A}+\tilde{A}(t)\\
 C(t)&=&\bar{C}
\end{array}\end{equation}
where
\begin{equation*}
\tilde{A}(t)=\frac{\partial \bar{f}}{\partial x}(\hat{x}(t),u(t)),
\end{equation*}
and
\begin{equation*}
\frac{\partial \bar{f}_{ir}}{\partial x_{jk}}(\hat{x}(t))=0
\end{equation*}
if $k>r$.

First we show that there exists $m_1>0$ such that for all $t\ge 0$
\begin{equation}\begin{array}{lllllllllllllllllllllllllllllllll}
P(t)&\le& m_1I^{n \times n}.
\end{array}\end{equation}
Consider the optimal control problem of minimizing
\begin{eqnarray*}
\xi'(0)P(0)\xi(0)+ \int_0^t \xi'(s)\Gamma\xi(s)+ \mu'(s)\mu(s)\ ds
\end{eqnarray*}
subject to
\begin{eqnarray*}
\dot{\xi}(s)&=& -A'(s) \xi(s) -\bar{C}'\mu(s)\\
\xi(t)&=& \zeta.
\end{eqnarray*}
It is well-known that the optimal cost is 
\begin{eqnarray*}
\zeta'P(t)\zeta
\end{eqnarray*}
where $P(t)$ is the solution of (\ref{ekf}).

Following \cite{GHO92} for $\theta>0$ we define $S(\theta)$
as the  solution of
\begin{eqnarray} \label{St}
\bar{A}'S(\theta) +S(\theta) \bar{A}-\bar{C}'\bar{C}&=&-\theta S(\theta) .
\end{eqnarray}
It is not hard to see that $S(\theta)$ is positive definite for $\theta>0$
as  it satisfies the Lyapunov equation
\begin{eqnarray*}
\left(-\bar{A} -{\theta \over 2}I\right)'S(\theta) 
+S(\theta)\left(-\bar{A} - {\theta \over 2}I\right)&=&  -\bar{C}'\bar{C}
\end{eqnarray*}
where $\bar{C},\left(-\bar{A} - {\theta \over 2}I\right)$ is an observable  pair and $\left(-\bar{A} - {\theta \over
2}I\right)$ has all eigenvalues equal to $-{\theta \over 2}$.
It follows from (\ref{St}) that
\begin{eqnarray*}
S_{ij,\rho\sigma}(\theta)&=&\frac{S_{ij,\rho\sigma}(1)}{\theta^{j+\sigma-1}}
=\frac{(-1)^{j+\sigma}}{\theta^{j+\sigma-1}}
\left(\begin{array}{c} j+\sigma-2 \\j-1 \end{array} \right).
\end{eqnarray*}

Let $T(\theta)=S^{-1}(\theta)>0$ then
\begin{eqnarray*}
T_{ij,\rho\sigma}(\theta)&=&{\theta^{j+\sigma-1}}{T_{ij,\rho\sigma}(1)}
\end{eqnarray*}
and $T(\theta)$ satisfies the Riccati equation
\begin{eqnarray*}
-\bar{A}T(\theta) -T(\theta) \bar{A}'+T(\theta)\bar{C}'\bar{C}T(\theta)&=&\theta T(\theta) .
\end{eqnarray*}

We apply the suboptimal control $\mu=-\bar{C} T(\theta) \xi$ to the above optimal control problem
and conclude that 
\begin{eqnarray} \label{subopt}
\zeta'P(t)\zeta & \le & \xi'(0)P(0)\xi(0)+ \int_0^t \xi'(s)\left(\Gamma+T(\theta)\bar{C}' \bar{C}
T(\theta)\right)\xi(s)\ ds
\end{eqnarray}
where
\begin{eqnarray*}
\dot{\xi}(s)&=& \left(-A'(s)  +\bar{C}'\bar{C} T(\theta) \right) \xi(s)\\
\xi(t)&=& \zeta.
\end{eqnarray*}

Now
\begin{eqnarray*}
\frac{d}{ds} \xi'(s)T(\theta) \xi(s)&=& 
\xi'(s)\left( \theta T(\theta)+ T(\theta)\bar{C}'\bar{C}T(\theta) \right)\xi(s)\\
&&
-\xi'(s)\left( \tilde{A}(s)T(\theta)+T(\theta)\tilde{A}'(s) \right)\xi(s).
\end{eqnarray*}
Because of the Lipschitz condition (\ref{LipGr}) we conclude that
$$ |A(s)| \le L$$
and 
$$ |\tilde{A}(s)| \le L+ |\bar{A}|. $$
From the form of $\tilde{A}(s)$ and $T(\theta)$ we conclude that
\begin{eqnarray*}
\left(\tilde{A}(s)T(\theta)\right)_{ij,\rho\sigma}&=&
 O(\theta)^{j+\sigma-1}
\end{eqnarray*}
while on the other hand
\begin{eqnarray*}
\theta T_{ij,\rho\sigma}(\theta)&=&{\theta^{j+\sigma}}{T_{ij,\rho\sigma}(1)}.
\end{eqnarray*}
Hence we conclude that for any $\alpha>0$ there exists $\theta$ sufficiently large
so that
\begin{eqnarray*}
\theta T(\theta)+ T(\theta)\bar{C}'\bar{C}T(\theta) -
 \tilde{A}(s)T(\theta)-T(\theta)\tilde{A}'(s) \ge \alpha I^{n\times n}.
\end{eqnarray*}
Therefore for $0\le s \le t$
\begin{eqnarray*}
\xi'(s)T(\theta) \xi(s)&\le& e^{\alpha(s-t)}\zeta'\zeta
\end{eqnarray*}
Now there exists $m_2(\theta)>0 $ such that
$$\xi'(s) \xi(s)\le m_2(\theta)\xi'(s)T(\theta) \xi(s)$$
so we conclude that
$$\xi'(s) \xi(s)\le m_2(\theta)e^{\alpha(s-t)}\zeta'\zeta.$$

There exist constants $m_3>0,m_4(\theta)>0$ such that
\begin{eqnarray*} 
P(0)&\le& m_3 I^{n \times n}\\
\Gamma+T(\theta)\bar{C}' \bar{C}
T(\theta)&\le& m_4(\theta) I^{n \times n}
\end{eqnarray*}
From (\ref{subopt}) we obtain the desired conclusion,
\begin{eqnarray*} \label{s1}
\zeta'P(t)\zeta & \le & m_3 e^{-\alpha t}\zeta'\zeta+m_4(\theta)\int_0^t e^{\alpha(s-t)}\zeta'\zeta\ ds\\
\zeta'P(t)\zeta & \le & m_3 \zeta'\zeta+m_4(\theta)\int_{-\infty}^t e^{\alpha(s-t)}\zeta'\zeta\ ds\\
\zeta'P(t)\zeta & \le &
{m_3  +m_4(\theta)\over \alpha}\zeta'\zeta
\end{eqnarray*}

Define
$$ Q(t)= P^{-1}(t)$$
then $Q$ satisfies
\begin{equation}\begin{array}{lllllllllllllllllllllllllllllllll}
\dot{Q}(t)&=&-A'(t)Q(t)-
Q(t)A(t)-Q(t)\Gamma Q(t) +\bar{C}'\bar{C}\\
Q(0)&=& P^{-1}(0) >0
\end{array}\end{equation}
Next we show that there exists $m_5>0$ such that for all $t\ge 0$
\begin{eqnarray*}
Q(t)&\le& m_5 I^{n \times n}.
\end{eqnarray*}
This will imply that 
\begin{eqnarray} \label{Pbelow}
P(t)&\ge& {1\over m_5} I^{n \times n}.
\end{eqnarray}

Consider the optimal control problem of minimizing
\begin{eqnarray*}
\xi'(0)Q(0)\xi(0)+ \int_0^t \xi'(s)\bar{C}'\bar{C}\xi(s)+ \mu'(s)\mu(s)\ ds
\end{eqnarray*}
subject to
\begin{eqnarray*}
\dot{\xi}(s)&=& A(s) \xi(s) +G\mu(s)\\
\xi(t)&=& \zeta.
\end{eqnarray*}
It is well-known that the optimal cost is 
\begin{eqnarray*}
\zeta'Q(t)\zeta
\end{eqnarray*}
where $Q(t)$ is the solution of (\ref{ekf}).

We use the suboptimal control 
$$ \mu(s)=  G^{-1}(\alpha I^{n\times n}-A(s))\xi(s)$$
so that the closed loop dynamics is
\begin{eqnarray*}
\dot{\xi}(s)&=&\alpha \xi(s)\\
\xi(t)&=&\zeta\\
\xi(s)&=&e^{\alpha(s-t)}\zeta.
\end{eqnarray*}
From this we obtain the desired bound
\begin{eqnarray*}
\zeta'Q(t)\zeta&\le& \xi'(0)Q(0)\xi(0)\\&&+ 
\int_0^t \xi'(s)\left(\bar{C}'\bar{C}+(\alpha I^{n\times n}-A'(s))\Gamma^{-1}(\alpha I^{n\times n}-A(s)) 
   \right)\xi(s) \ ds \\
\zeta'Q(t)\zeta&\le& e^{-2\alpha t}\zeta'Q(0)\zeta+ \int_0^t e^{2\alpha(s-t)}\zeta'\left(\bar{C}'\bar{C}+(\alpha+L)^2
\Gamma^{-1}    \right)\zeta \ ds\\
\zeta'Q(t)\zeta&\le& \left(m_6+{m_7\over 2\alpha}\right)\zeta'\zeta
\end{eqnarray*}
where
\begin{eqnarray*}
Q(0)&\le&m_6 I^{n\times n}\\
\bar{C}'\bar{C}+(\alpha+L)^2
\Gamma^{-1} &\le&m_6 I^{n\times n}.
\end{eqnarray*}

Now let $x(t), u(t), y(t)$ be a trajectory of the system (\ref{nl3}) starting at $x^0$.  Let  $\hat{x}(t)$
be the trajectory of the extended Kalman filter (\ref{ekf}) starting at $\hat{x}^0$ and
$\tilde{x}(t)=x(t) -\hat{x}(t)$, $\tilde{x}^0=x^0 -\hat{x}^0$.
Then
\begin{eqnarray*}
\frac{d}{dt}\tilde{x}'(t)Q(t)\tilde{x}(t)&=&-\tilde{x}'(t)\left(\bar{C}'\bar{C} + Q(t)\Gamma Q(t)\right)\tilde{x}(t)\\
&&+ 2\tilde{x}'(t)Q(t)\left(\bar{f}(x(t),u(t))-\bar{f}(\hat{x}(t),u(t))-\tilde{A}(t)\tilde{x}(t)\right).
\end{eqnarray*}
Now following \cite{BBJ88}
\begin{eqnarray*}
\bar{f}(x(t),u(t))-\bar{f}(\hat{x}(t),u(t))-\tilde{A}(t)\tilde{x}(t)&=&\int_0^1\int_0^1 r \frac{\partial^2 \bar{f}}{\partial
x_i
\partial x_j}(\hat{x}(t)+rs\tilde{x}(t),u(t))\tilde{x}_i(t)\tilde{x}_j(t)\ ds \ dr\\
&\le& L |\tilde{x}(t)|^2.
\end{eqnarray*}
Since $G$ is invertible there exists $m_7>0$ such that
\begin{eqnarray*}
\Gamma&\ge& m_7 I^{n \times n}
\end{eqnarray*}
and so 
\begin{eqnarray*}
\frac{d}{dt}\tilde{x}'(t)Q(t)\tilde{x}(t)&\le &-{ m_7 \over m_1^2}|\tilde{x}(t)|^2+m_5 L
|\tilde{x}(t)|^3\\
&\le &-{ m_7 \over m_1^2m_5}\tilde{x}'(t)Q(t)\tilde{x}(t)+{m_5 L }
(m_1\tilde{x}'(t)Q(t)\tilde{x}(t))^{3\over 2}
\end{eqnarray*}
If 
$$ (\tilde{x}'(t)Q(t)\tilde{x}(t))^{1\over 2}< { m_7 \over 2 m_1^{7\over 2} m_5^2 L}$$
then
\begin{eqnarray*}
\frac{d}{dt}\tilde{x}'(t)Q(t)\tilde{x}(t)
&\le&-{ m_7 \over 2  m_1^2 m_5} \tilde{x}'(t)Q(t)\tilde{x}(t)
\end{eqnarray*} 
so $\tilde{x}'(t)Q(t)\tilde{x}(t) \to 0$ exponentially as $t \to \infty$.  Therefore if 
$$ (\tilde{x}'(t)Q(t)\tilde{x}(0))^{1\over 2}< { m_7 \over 2 m_1^{7\over 2} m_5^2 L}$$
the extended Kalman filter converges.

\section{Conclusions}
\setcounter{equation}{0}

The above result does not follow from that of
Baras, Bensoussan and James \cite{BBJ88}.  To show local convergence of  the extended Kalman
filter they required "uniform detectability".  They define this as follows.
The system
\begin{equation}\begin{array}{lllllllllllllllllllllllllllllllll} \label{nlbbj}
\dot{x}&=& f(x)\\
y&=& Cx
\end{array}\end{equation} 
is uniformly detectable if there exists a bounded Borel matrix-valued function
$\Lambda(x)$ and a constant $\alpha>0$ such that for all $x,\xi\in \Real^n$
\begin{eqnarray*}
\xi'\left(\frac{\partial f}{\partial x}(x)+\Lambda(x)C\right) \xi &\le&
-\alpha |\xi|^2.
\end{eqnarray*}

This is a fairly restrictive condition as not all observable linear
systems are uniformly detectable.  Consider
\begin{eqnarray*}
\dot{x}&=&  \left[ \begin{array}{ccccccccc} 0&1\\a_1&a_2\end{array}\right] x\\
y&=&\left[ \begin{array}{ccccccccc} 1&0\end{array}\right] x.
\end{eqnarray*}
Suppose
$$ \Lambda(x)= \left[ \begin{array}{ccccccccc} \lambda_1(x)\\ \lambda_2(x) \end{array}\right]$$
then
$$\frac{\partial f}{\partial x}(x)+\Lambda(x)C=\left[ \begin{array}{ccccccccc} \lambda_1(x)&1\\\lambda_2(x)+a_1&a_2\end{array}\right].$$
If $a_2>0$ and $\xi' = \left[ \begin{array}{ccccccccc} 0&1 \end{array}\right]$ then
\begin{eqnarray*}
\xi'\left(\frac{\partial f}{\partial x}(x)+\Lambda(x)C\right) \xi &=& a_2 > 0
\end{eqnarray*}
so the system is not uniformly detectable.  This system does satisfies the conditions
of Theorem \ref{mt} so  an extended Kalman filter would converge locally.  Since the system is linear,
an extended Kalman filter is also a Kalman filter that converges globally 

An example \cite{KD96} of a highly nonlinear problem where an EKF may fail to 
converge is
\begin{equation}\begin{array}{lllllllllllllllllllllllllllllllll}
\dot x &=&f(x)= x(1-x^2) \\
y&=& h(x)=x^2-  x / 2
\end{array}\end{equation}
where $x,y \in \Real$.
The system is observable as $y, \dot{y}, \ddot{y}$ separate points but it is not uniformly observable.
  The
dynamics has two stable equilibria at $x=\pm1$ and an unstable equilibrium
at $x=0$.  Under certain initial conditions, the extended Kalman filter  fails to 
converge.  Suppose the $x^0=1$ so $x(t)=1$  and $y(t)=1/2$ for
all $t\geq 0$.  But $h(-1/2)=1/2$ so if  $\hat{x}^0\leq -1/2$ the extended
Kalman filter will not converge. To see this notice that when $\hat{x}(t)= -1/2$, the
term $y(t)-h(\hat{x}(t) )=0$ so $\dot{\hat{x}}=f(\hat{x}(t) )=f(-1/2)=-3/8$.
Therefore $\hat{x}(t)\leq -1/2$ for all $t\geq 0$.

\section{Dedication}
\setcounter{equation}{0}
  This paper is dedicated to my esteemed colleague and good friend, Professor Anders Lindquist
on the occasion of his $60^{th}$ birthday.

\end{document}